\documentclass{article}
\usepackage{graphicx}
\def\Re{\hbox{Re\kern 1pt}}
\def\ch{1}
\def\fem{2}
\def\serkh{3}
\def\wasow{4}
\def\newman{5}
\def\stahl{6}
\def\atap{7}
\def\aaa{8}
\def\gopalt{9}
\def\fhm{10}
\def\NA{{11}}
\def\fenics{{12}}
\def\firedrake{{13}}
\def\ifiss{{14}}
\def\pltmg{{15}}
\def\qbx{{16}}
\def\bb{{17}}
\def\fair{{18}}
\def\liu{{19}}
\def\rokb{{20}}
\def\bbb{{21}}
\def\evans{{22}}

\begin{document}

\noindent
{\Large\bf New Laplace and Helmholtz solvers}

\bigskip

\noindent
Abinand Gopal and Lloyd N. Trefethen\hfill\break
Mathematical Institute, University of Oxford, Oxford OX2 6GG, UK\hfill\break
{\tt gopal@maths.ox.ac.uk} \kern 1.5pt and \kern .5pt{\tt trefethen@maths.ox.ac.uk}

\vskip .4in

\noindent\textsf{\textbf{{\em Abstract.}  New numerical algorithms
based on rational functions are introduced that can solve certain
Laplace and Helmholtz problems on two-dimensional domains with
corners faster and more accurately than the standard methods
of finite elements and integral equations.  The new algorithms
point to a reconsideration of the assumptions underlying existing
numerical analysis for partial differential equations.}}

\vskip .4in

\noindent{\bf 1. Laplace equation}

The Laplace and Helmholtz equations are the basic partial
differential equations (PDE\kern .3pt s) of potential theory and
acoustics, respectively [\ch].  Suppose a region $\Omega$ bounded
by a polygon $P$ is given and (to begin with the Laplace case)
we seek the unique function $u(x,y)$ that satisfies $\Delta
u = \partial^2u/\partial x^2 + \partial^2u/\partial y^2 = 0$
in $\Omega$ and matches a given function $h$ at the boundary.
More generally, $P$ might have curved sides meeting at corners,
and the boundary data might involve derivatives as well as
function values.  It is convenient to represent the coordinates
by a complex variable $z=x+iy$, so we write the boundary condition
as $u(z) = h(z)$ for $z\in P$.

The standard techniques for solving such a problem numerically
are the finite element method (FEM) [\fem] and boundary integral
equations [\serkh].  Yet these methods face a challenge in
calculating accurate solutions because of singularities at the
corners [\wasow].  The mathematical basis of a new algorithm for
meeting this challenge is a theorem in the field of approximation
theory published by Donald Newman in 1964 [\newman].  Newman
considered the problem of approximation of $f(x) = |x|$ on the
interval $[-1,1]$ by a rational function, that is, a quotient of
polynomials $r(z) = p(z)/q(z)$.  He showed that whereas polynomial
approximations converge at best at the very slow rate $\|f-p_n\|
= O(n^{-1})$, rational approximations can achieve much faster
``root-exponential'' convergence $\|f-r_n\| = O(\exp(-C\sqrt n\kern
1pt ))$ with $C>0$ [\stahl,\atap].  Here $n$ is the degree of a
polynomial or rational function, which is defined in the latter
case as the maximum of the degrees of $p$ and $q$.

Our new algorithm achieves root-exponential
convergence for solving PDE\kern .3pt s by approximating $u(z)$ by
the real part of a rational
function, $u(z) = \Re r(z)$.  Any such approximation is a harmonic
function in $\Omega$, i.e., a solution of the Laplace equation,
provided $r$ has no poles in $\Omega$. 
Finding rational approximations to given data is a difficult
nonlinear problem in general [\aaa].\ \ Here, however,
we know that the dominant singularities of the function
to be approximated are located at the vertices of $P$.
This suggests the idea, motivated by Newman's result and
related computational experience [\aaa], of prescribing poles of $r$ outside
$\Omega$ a priori in a configuration with exponential
clustering near each vertex.  
Specifically, our
rational functions are represented in the partial fraction form
\begin{equation}
r(z) = \sum_{j=1}^{N_1}{a_j\over z-z_j} + \sum_{j=0}^{N_2}b_j z^j,
\label{partfrac}
\end{equation}
with the poles $\{z_j\}$ fixed and the real and imaginary parts of
$\{a_j\}$ and $\{b_j\}$ as unknowns ($N=2N_1 + 2N_2+1$ degrees of
freedom in total, since $b_0$ can be taken to be real).  One might
expect that the pole locations would have to be delicately
chosen to be effective.  However, in a new mathematical result
to be published elsewhere, we have proved that a straightforward
realization of this idea, relying on no unknown parameters, is
enough to guarantee the existence of approximations of the form (1)
with root-exponential convergence.  To find such approximations
computationally, one exploits the fact that since the poles are
prescribed, the problem is linear.  Expansion coefficients are
found by least-squares fitting in sample points on the boundary, a
routine problem of linear algebra involving a matrix of dimensions
about $3N\times N$.\ \ Starting from a small value such as $50$,
$N$ is increased systematically until a prescribed accuracy
is achieved.  In typical problems on polygons with up to $8$
vertices, we find that $N\approx 1000$ suffices to give accuracy
to 8 digits, complete with an accuracy guarantee derived from the
maximum principle, which asserts that the error at each point in
the interior is bounded by the maximal error on the boundary.

\begin{figure}[h]
\begin{center}
\bigskip\medskip
\includegraphics[scale=1.15]{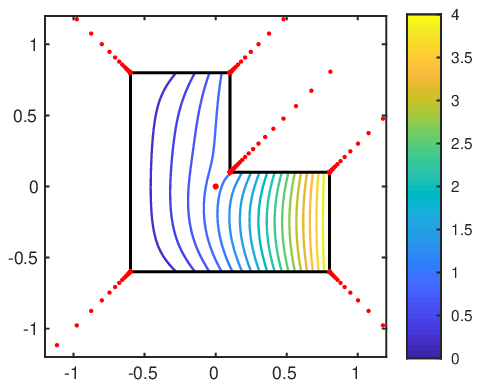}
\end{center}
{\small{\bf Fig.~1.  Laplace equation in an L-shaped region.}
Poles are clustered exponentially near each vertex, and then
a least-squares problem is solved on the boundary $P$ to find
expansion coefficients for a global representation
$(\ref{partfrac})$ of a solution accurate to 10 digits.
The dot in the interior
of the domain marks the expansion point of
the polynomial part of $(\ref{partfrac})$.\par}
\end{figure}

Figure 1 shows the solution computed by the new Laplace equation
algorithm in an L-shaped region, a famous test geometry for
problems of this kind [\fhm].  For a problem that is generic in
the sense of having singularities at the corners, the simplest
choice of boundary data is $u(z) = h(z) = [\Re z]^2 = x^2$.
Our code approximates the solution successively with $N=42$,
$82,$ $138,\dots,$ $1002$ degrees of freedom, at which point
10-digit accuracy is achieved.  Figure 2 shows the accuracy as a
function of $\sqrt N$, revealing a straight line corresponding
to root-exponential convergence.  For a wide range of Laplace
problems on polygons, this performance is representative.
Typically we solve a problem in ${<}\kern 1pt 1$\kern .3pt s on
a laptop running MATLAB in 16-digit floating-point arithmetic,
and then each evaluation of the solution, with guaranteed accuracy
all the way to the corners, takes a few tens of $\mu$s.  For the
problem of Fig.~1, the evaluation of $u$ to 10-digit accuracy at
$10^4$ points in $\Omega$ takes 0.3\kern .3pt s.

\begin{figure}
\begin{center}
\includegraphics[scale=1.0]{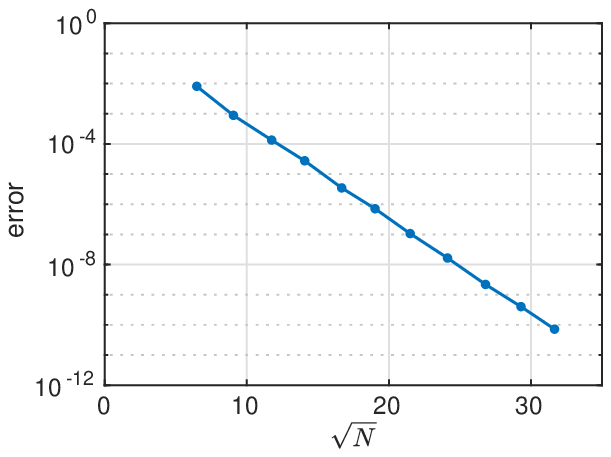}
\end{center}
{\small{\bf Fig.~2.  Convergence.}  The
maximal error as a function of number of degrees
of freedom $N$ for the problem of Fig.~1 shows
root-exponential convergence.  Similar convergence is observed
for other domains with corners.\par}
\end{figure}

\bigskip

\noindent{\bf 2. Comparison with other methods}

The numerical solution of PDE\kern .3pt s has been at the heart of
scientific computing since the 1950\kern .3pt s, and the Laplace
equation in a planar domain is as fundamental a problem in this area
as any.  There are two main classes of methods for solving such
problems: the finite element method (FEM), and boundary integral
equations.  To gather information on how the new approach compares
with these, in November 2018 we posed the L-domain problem just
presented as a challenge to the international numerical analysis
community via the email list NA Digest [\NA].   Specifically, we
asked for a computation of $u(0.99,0.99)$ to 8 digits of accuracy
(the exact value is $1.02679192610\dots$).  This led to responses
from about twenty experts around the world.  As expected, about
half the responses recommended FEM algorithms and software such as
FEniCS [\fenics], Firedrake [\firedrake], IFISS [\ifiss], and PLTMG
[\pltmg].  In this approach, the set of functions on $\Omega$ is
approximated by a finite dimensional linear space, and then one
solves a matrix problem to find a good candidate in that space.
Typically the matrices are large, sparse, and ill-conditioned.
The FEM is noted for its great flexibility, enabling problems much
more complicated than ours to be solved effectively, in three
as well as two dimensions.  However, it is not very efficient
near singularities, and whereas all our respondents were able
to calculate a solution to 2--4 digits of accuracy, only one
came close to 8 digits.  For example, one researcher used 158,997
5th-order triangular finite elements near the reentrant corner and
achieved 6 correct digits.  Our assessment is that it is possible
to solve the Laplace equation on a region with corners to high
accuracy by the FEM, but this entails a significant computation
requiring a high level of expertise and tools.

The other well-known methods are boundary integral equations, which
are advantageous because of good conditioning and because
the solution is represented in one dimension, along the
boundary, rather than two dimensions in the domain [\serkh].
Here one begins by solving an integral equation to determine
a density function $\rho(z)$, such as this
Fredholm equation of the second kind:
\begin{equation}
h(z) = - \pi \rho(z) + 
\int_P K(z,\zeta) \rho(\zeta) |d\zeta|, \quad
	K(z,\zeta) = {\cos(\angle(z-\zeta,\nu(\zeta)))\over |z-\zeta|},
\end{equation}
where $\angle(z-\zeta,\nu(\zeta))$ is the angle between $z-\zeta$ and
the inward normal to $P$ at $\zeta$.
This is a ``double-layer potential'' formulation,
in the standard terminology, whose solution $\rho(z)$ can
be interpreted as a distribution of dipole charge density
along the boundary.  Once $\rho$
is found, the solution to the Laplace
problem is evaluated at a point $z$ by computing an integral:
\begin{equation}
u(z) = \int_P K(z,\zeta) \rho(\zeta) |d\zeta|.
\end{equation}
In general $\rho$ will be smooth along the sides and singular at
the corners, introducing challenges in evaluating the integrals.
Additionally, the kernel $K$ is singular, so that even away
from the corners, an accurate evaluation may be a nontrivial
task.  However, experts have developed powerful techniques for
such quadratures [\qbx], and four respondents to our inquiry
produced results accurate to the specified 8 digits or more.
We conclude that if one wants an accurate solution to Laplace
problems in corners, integral equations are the most effective
of the existing technologies.  On the other hand, there is not
much software available, and the solutions communicated to us
were produced by experts running their own codes.

Whereas integral equations make use of a continuous distribution
of dipoles $\rho(z)$ on the boundary, the new rational functions
method can be interpreted as using a finite sum of dipoles (delta
functions) beyond the boundary.  This has the advantage that
evaluation of $u(z)$ takes exactly the same form (1) regardless of
$z$, requiring no special calculations when $z$ is near a boundary
arc or a corner.  One may regard this as a zero-dimensional
representation of the solution, whose complexity is determined
only by the complexity of the function being represented, not
by that of the region $\Omega$ or its boundary $P$.  The most
closely related existing method goes by the name of the method
of fundamental solutions (MFS), also called the charge simulation
method, in which solutions are also approximated via finite sums
[\bb,\fair].  The MFS differs from our approach in that each
term is normally a monopole (point charge) rather than a dipole,
and it is not normally applied with exponential clustering to
achieve root-exponential convergence [\liu].

\bigskip

\noindent{\bf 3. Helmholtz equation}

The most important generalization of the Laplace equation is 
the Helmholtz equation $\Delta u + k^2 u = 0$,
which models time-harmonic propagation of acoustic or
electromagnetic waves at frequency $k$ [\rokb].\ \ For a Helmholtz
problem in a domain exterior to a scattering body
with wave fields satisfying the Sommerfeld radiation
condition, we modify (\ref{partfrac}) to
\begin{equation}
\sum_{j=1}^{N_1} a_j H_0(k|z-z_j|) +
\sum_{j=1}^{N_1} b_j H_1(k|z-z_j|){z-z_j\over |z-z_j|} +
\sum_{j=0}^{N_2} c_j H_j(k|z|){z^j\over |z|^j},
\label{partfrac2}
\end{equation}
where $k$ is fixed and $H_j$ are Hankel functions of the first
kind.  We cluster singularities exponentially near corners in
the interior of the scatterer.  Interesting problems now require
complex expansion coefficients to cancel a signal incident at
the boundary; the ``sound soft'' case has zero total field at the
boundary, and the ``sound hard'' case has zero normal derivative.
Figure 3 shows two sound soft example problems solved to 4-digit
accuracy [\bbb].  In each case the solution took about 2\kern .3pt
s on a laptop, with about 250 $\mu$s for each point evaluation
afterwards.
The convergence rate is again root-exponential, as shown in Fig.~4.

\begin{figure}[t]
\begin{center}
\includegraphics[scale=.64]{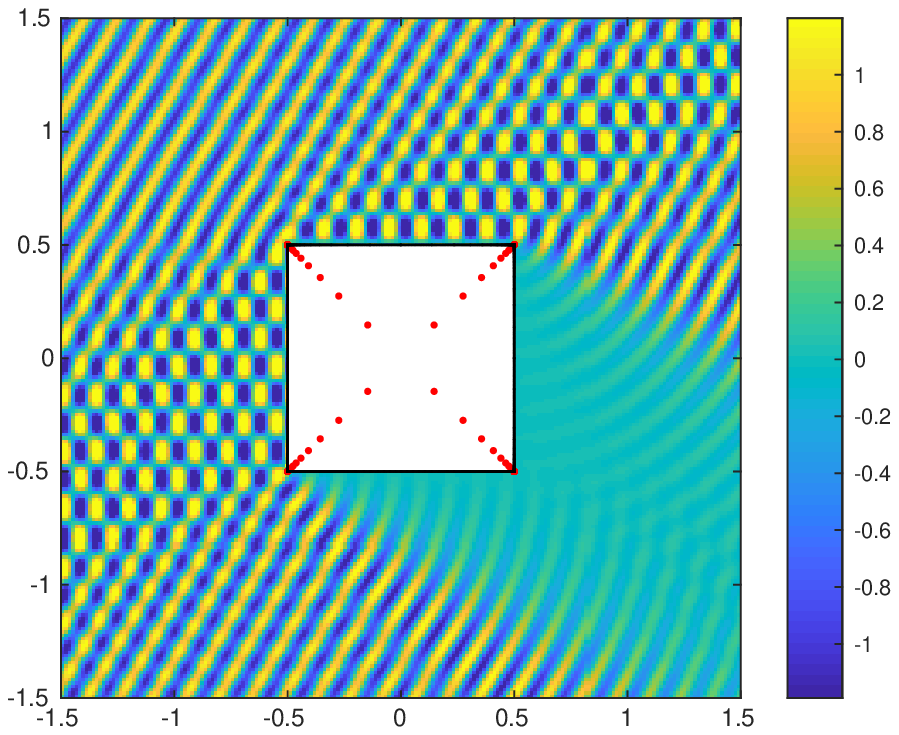}
~~~~~~~~\includegraphics[scale=.64]{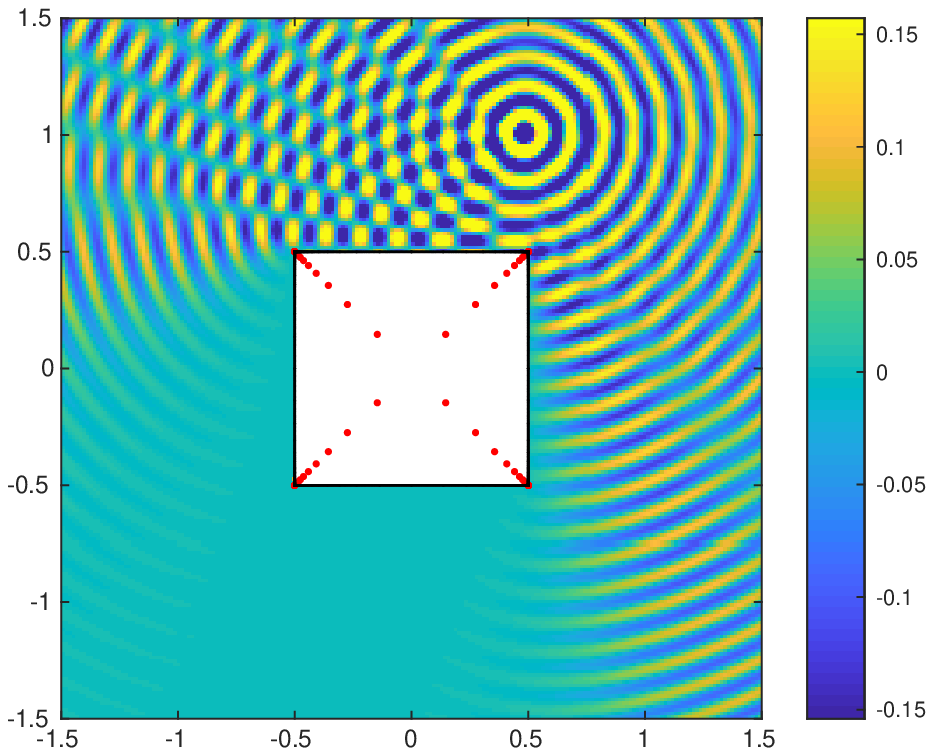}
\end{center}
{\small{\bf Fig.~3.  Helmholtz equation.}
Analogous computations to Fig.~1 for the Helmholtz
equation $\Delta u + k^2 u = 0$ with $k=50$ in the exterior of a square.
The incident signal on the left  is a plane wave oriented
at $30^\circ$, and on the right, a point
	oscillation $H_0(k|z-z_0|)$ situated at $z_0 = 1/2+i$.\par}
\end{figure}

\begin{figure}[t]
\begin{center}
\bigskip\medskip
\includegraphics[scale=1.0]{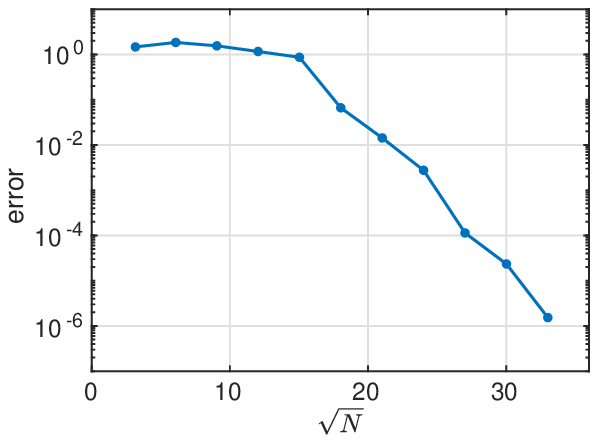}
\end{center}
{\small{\bf Fig.~4.  Convergence.}  Root-exponential convergence for
the Helmholtz equation, after an initial $k$-dependent transient needed
to resolve the wave.
The quantity measured is the maximum error on the boundary for
the first example of Fig.~3. \par}
\end{figure}

\bigskip

\noindent{\bf 4. Discussion}

Many issues will need to be addressed as the new methods
are developed.  One is the effective choice of exponential
clustering parameters and numbers of singularities near each
vertex, which does not matter in a certain theoretical sense
but may be important in practice since poor choices lead to
numerical instability.  Another is the treatment of slits and
other geometric complications, including the treatment of highly
elongated domains, where it may be important to generalize the
$N_2$ sums of (\ref{partfrac}) and (\ref{partfrac2}) to several
sums centered at different points.  For Helmholtz problems with
high frequency ($k\gg 1$), a simple $N_2$ sum may cease to
be competitive for the smooth part of the problem with other
established methods.  And of course there is the question
of extension to the three space dimensions, where again,
existing methods do not always well exploit the smoothness
of most problems away from certain singular points and curves.
For the Laplace equation, this entails a move from the logarithmic
potential associated with rational functions and the plane to the
inverse-linear potential of Newtonian mechanics.  In principle,
methods with fractional-exponential convergence analogous to the
methods presented here should exist in 3D, but we have no view as
to whether they will prove competitive in practice.

There is a historical context that may shed light on the lack of
previous literature on solving PDE\kern .3pt s by means of rational
functions and their generalizations.  The mathematicians of the
nineteenth century focused on functions that were analytic or
piecewise analytic, i.e., smooth and representable by convergent
Taylor series.  Most physical applications are of this kind.
In the twentieth century, however, mathematicians turned to new
challenges of less smooth functions, developing advanced tools
for analyzing fine distinctions of regularity (i.e., smoothness).
Overwhelmingly, such tools became the standard framework for
computational mathematicians too, and in particular, FEM experts
almost invariably derive and analyze their algorithms in the
language of Sobolev spaces [\evans], in which precise distinctions
are made, for example, between a function with one derivative
of smoothness and a function with one-and-a-half derivatives.
This kind of numerical analysis brings with it a bias towards
low-accuracy methods tuned to problems with limited smoothness.
The simplicity and speed of the new methods proposed here are a
reminder that there is also a place for numerical analysis based
on less pessimistic smoothness assumptions.

This note has introduced a new class of numerical methods for
the fast and accurate solution of certain PDE\kern .3 pt s.
The method captures singularities to high accuracy without having
to analyze them and delivers a global representation of the solution
by a single formula.  To assess its prospects, one may divide the
scale of PDE problems into small, medium, and large.  For small
or ``toy'' problems such as illustrated here, the new method
appears to be faster than existing methods for high-accuracy solutions
(e.g., 5--10 digits).  At the other extreme, for large problems
of computational engineering, often in complex geometries, the
new approach would probably be unworkable (and in any case would
require extension to three dimensions).  It is the middle range,
which one might associate with computational science more than
computational engineering, that will be most important in the
long run.  We believe the new approach holds promise for such
problems and that developing it is a challenge for the years ahead.
The established methods of finite elements and integral equations,
of course, have benefited from half a century's head start.

\bigskip\bigskip

\noindent{\bf Acknowledgements.} This work was initiated during
September 2018 while the authors were visiting New York University.
We thank Alex Barnett, Timo Betcke, Leslie Greengard, Dave
Hewett, Yuji Nakatsukasa, Vladimir Rokhlin, Kirill Serkh, and
Euan Spence for their advice.

\newpage
\noindent{{\bf References}}

\medskip

{\parindent=0pt

{
[\ch]
R. Courant and D. Hilbert, {\em Methods
of Mathematical Physics}
(Interscience, 1966).

[\fem]
D. Braess, {\em Finite Elements: Theory, Fast Solvers,
and Applications in Solid Mechanics} (Cambridge U. Press, 2007).

[\serkh]
K. Serkh and V. Rokhlin,
On the solution of elliptic partial
differential equations on regions with corners,
{\em J. Comput.\ Phys.}\ {\bf 305}, 150--171 (2016).

[\wasow]
W. Wasow,
Asymptotic development of
the solution of Dirichlet's problem at analytic corners,
{\em Duke Math.\ J.} {\bf 24}, 47--56 (1957).

[\newman]
D. J. Newman,
Rational approximation to $|x|$,
{\em Mich.\ Math.\ J.} {\bf 11}, 11--14 (1964).

[\stahl]
H. Stahl,
Best uniform rational approximation of $|x|$ on $[-1,1]$,
{\em Russian Acad.\ Sci.\ Sb. Math.} {\bf 76}, 461--487 (1993).

[\atap]
L. N. Trefethen, {\em Approximation Theory and Approximation Practice}
(SIAM, 2013).

[\aaa]
Y. Nakatsukasa, O. S\`ete, and L. N. Trefethen,
The AAA algorithm for rational approximation,
{\em SIAM J. Sci.\ Comput.}\ {\bf 40}, A1494--A1522 (2018).

[\gopalt]
A. Gopal and L. N. Trefethen,
Representation of conformal maps by rational functions,
{\em Numer.\ Math.,} to appear.

[\fhm]
L. Fox, P. Henrici, and C. Moler,
Approximations and bounds for eigenvalues of elliptic operators,
{\em SIAM J. Numer.\ Anal.}\ {\bf 4}, 89--102 (1967).

[\NA]
L. N. Trefethen, 8-digit Laplace solutions
on polygons?\ \ Posting on NA Digest at\hfill\break
{\tt www.netlib.org/na-digest-html} (29 November 2018).

[\fenics]
{\tt https://fenicsproject.org}.

[\firedrake]
{\tt https://www.firedrakeproject.org}.

[\ifiss]
H. C. Elman, A. Ramage, and D. J. Silvester,
IFISS: A computational laboratory for investigating
incompressible flow problems,
{\em SIAM Rev.} {\bf 56}, 261--273 (2014).

[\pltmg]
R. E. Bank, PLTMG 13.0, {\tt www.netlib.org} (2018).

[\qbx]
A. Kl\"ockner, A. Barnett, L. Greengard, and M. O'Neil,
Quadrature by expansion: A new method
for the evaluation of layer potentials,
{\em J. Comput.\ Phys.} {\bf 252}, 332--349 (2013).

[\bb]
A. H. Barnett and T. Betcke,
Stability and convergence of the method of fundamental solutions for
Helmholtz problems on analytic domains,
{\em J. Comput.\ Phys.} {\bf 227}, 7003--7026 (2008).

[\fair]
G. Fairweather and A. Karageorghis, 
The method of fundamental solutions for elliptic
boundary value problems,
{\em Adv.\ Comput.\ Math.} {\bf 9}, 69--95 (1998).

[\liu]
Y. Liu, {\em The Numerical Solution of
Frequency-Domain Acoustic and Electromagnetic Periodic Scattering
Problems,}  PhD thesis, Dartmouth College (2016).

[\rokb]
V. Rokhlin,
Rapid solution of integral equations of
scattering theory in two dimensions,
{\em J. Comput.\ Phys.} {\bf 86}, 414--439 (1990).

[\bbb]
A. H. Barnett and T. Betcke,
An exponentially convergent nonpolynomial
finite element method for time-harmonic scattering from polygons,
{\em SIAM J. Sci.\ Comput.} {\bf 32}, 1417--1441 (2010).

[\evans]
L. C. Evans, {\em Partial Differential Equations,}
Amer.\ Math.\ Soc., 1998.

\par}

\medskip
\smallskip

\end{document}